 \documentstyle[bezier,12pt]{article}
\def \R{{I\!\!\!\,R}}

\def\n0N{n=0,1,\cdots,N}
\def\n1N{n=1,\cdots,N}
\def\0nN{0\le n\le N}
\def\1nN{1\le n\le N}

\def\text#1{{\rm #1}}

\newtheorem{Theorem}{Theorem}[part]

\begin{document}
\date{February  19, 2008 (Revised); September 13, 2007 (Original)}
\title{\LARGE General Matrix-Valued Inhomogeneous Linear Stochastic Differential Equations
  and Applications}
\author{Jinqiao Duan \footnote{Supported by the NSF Grant 0620539.} \ \ \ \ \ Jia-an Yan
\footnote{Supported by the National Natural Science Foundation of
China (No. 10571167)and the National Basic Research Program of
China (973 Program) (No. 2007CB814902).  }\\
 \small Department of Applied Mathematics\  \ \ \ \ Institute of Applied
 Mathematics\\
\small \hskip40pt  Illinois Institute of Technology \hskip30pt  Academy of Math. and Systems Science\\
\small  \hskip40pt  Chicago, IL 60616, USA \hskip30pt   Chinese
Academy of
Sciences, Beijing 100080 \\
\small    Email:  \emph{duan@iit.edu}   \hskip30pt  Email:
\emph{jayan@amt.ac.cn}}

\small \maketitle {\bf Abstract: } The expressions of solutions
for general $n\times m$ matrix-valued inhomogeneous linear
stochastic differential equations are derived. This generalizes a
result of Jaschke (2003) for scalar inhomogeneous linear
stochastic differential equations. As an application, some $\R^n$
vector-valued inhomogeneous nonlinear stochastic differential
equations   are reduced to random differential equations,
facilitating pathwise study of the solutions.

\section{A Review of Stochastic Exponential Formulas } \setcounter{equation}{0}
\setcounter{Assumption}{0} \setcounter{Theorem}{0}
\setcounter{Proposition}{0} \setcounter{Corollary}{0}
\setcounter{Lemma}{0}

We first review some existing results about solution formulas for
linear stochastic differential equations (SDEs) or for their
integral formulations. Let $(\Omega,{\cal F},({\cal F} t),{\bf
P})$ be a standard stochastic basis. For the following stochastic
integral equation
$$
X_t=1+\int_0^t X_{s-} d Z_s, \eqno(1.1),
$$
where $Z$ is a semimartingale with $Z_0=0$, Dol\'eans-Dade (1970)
proved that the unique solution of (1.1) is given by
$$X_t=\exp \Big \{ Z_t-\frac12 \langle Z^c,
Z^c\rangle_ t \Big \} \prod_{0 <s \le t} (1+\Delta Z_s) e^{-\Delta
Z_ s}.\eqno(1.2)$$  In the literature  $X$ is called  {\it
Dol\'eans exponential (or stochastic  exponential) } of $Z$, and
is denoted  by ${\cal E}(Z)$.  The formula (1.2) is called the
{\it Dol\'eans (or stochastic) exponential formula}.

In an unpublished paper, Yoeurp and Yor (1977)  proved the
following result for the solution formula of scalar SDEs (see also
Revuz-Yor (1999) and Protter (2005) for the case where $Z$ is a
continuous semimartingale, and Melnikov-Shiryaev (1996)) for the
general case).

\begin{Theorem}  (\textbf{Yoeurp and Yor  1977} )  Let $Z$  and $H$ be semimartingales,  and $\Delta Z_s\ne -1$ for all $s\in
[0,\infty]$. Then the unique solution of the inhomogeneous scalar
linear SDE
$$X_t=H_t+\int_0^t X_{s-} d Z_s, \quad t \ge 0, \eqno(1.3)$$
is given by
$$X_t={\cal E} (Z)_t \Big \{ H_0 +\int_0^t {\cal E} (Z)_{s-}^{-1} dG_s\Big \},\eqno(1.4)$$
where
$$G_t=H_t-\langle H^c,Z^c\rangle_t-\sum_{0<s\le t}\frac{\Delta H_s\Delta
Z_s}{1+\Delta Z_s}.\eqno(1.5)$$
\end{Theorem}

Jaschke (2003) extended equation (1.3) to the case where  $H$ is an
adapted cadlag process, not necessarily a semimartingale. He proved
that in this case the solution of (1.3) is given by:
$$X_t=H_t-{\cal E} (Z)_t \int_0^t H_{s-}d({\cal E}
(Z)_s^{-1}).\eqno(1.6)$$ By using (1.6) Jaschke (2003) obtained a
new proof of (1.4).

On the other hand, Emery (1978)  considered the following $n\times
n$ matrix-valued stochastic equation
$$ U(t)=I+\int_0^t (d L(s))U(s-), \eqno(1.7),$$
where $I$ is an $n\times n$ identity, $L$ is a given $n\times n$
matrix-valued semimartingale with $L(0)=0$, such that $I+\Delta
L(s)$ is invertible for all $s\in [0,\infty]$, where $\Delta
L(s)=L(s)-L(s-)$. Emery proved that the equation (1.7) admits a
unique solution $U$, which is $n\times n$ matrix-valued
semimartingale. We call it the stochastic exponential of $L$ and
denote it  by ${\cal E}(L)$. However, there is no explicit
expression for such a stochastic exponential in general.

Jacod (1982) has studied the following inhomogeneous matrix-valued
stochastic integral equation
$$
X(t)=H(t)+\int_0^t (d L(s))X(s-), \eqno(1.8),
$$
where $L$ is a given $n\times n$ matrix-valued semimartingale with
$L 0=0$, such that $I+\Delta L(s)$ is invertible for all $s\in
[0,\infty]$, and $H$ is an $n\times m$ matrix valued semimartingale.
For an $n\times n$ matrix valued semimartingale $A$ and an $n\times
m$ matrix valued semimartingale $B$, we let
$$[A,B](t)=\langle A^c, B^c\rangle(t) + \sum_{0<s\le t}\Delta A(s)\Delta
B(s).$$ Here $A^c$ denotes its continuous martingale part defined
componentwise by $(L^c)^i_j=(L^i_j)^c$, and
$$\langle A^c, B^c\rangle^i_j=\sum_{k}\langle(A^c)^i k,(B^c)^k_j\rangle.\eqno(1.9)$$
Using these notations the result of Jacod (1982) implies the
following

\begin{Theorem} (\textbf{Jacod 1982}) The unique solution of (1.8) is given by
$$X(t)={\cal E} (L)(t) \Big \{ H(0) +\int_0^t {\cal E} (L)(s-)^{-1} dG(s)\Big \},\eqno(1.10)$$
where
$$G(t)=H(t)-\langle L^c,H^c\rangle(t)-\sum_{0<s\le t}({1+\Delta L(s)})^{-1}\Delta L(s)\Delta
H(s).\eqno(1.11)$$
\end{Theorem}

In particular, if $L$ and $H$ are continuous semimartingales, an
expression for the solution of (1.8) is given in Revuz-Yor (1999) as
follows:
$$X(t)={\cal E} (L)(t) \Big \{ H(0) +\int_0^t {\cal E} (L)(s)^{-1} (dH(s)-d[L,H](s))\Big\}.\eqno(1.12)$$

The objective of the present note is to generalize equation (1.8)
to the case where $H(t)$ is a given $n\times m$ matrix-valued
adapted cadlag process, not necessarily a semimartingale. We give
an expression of the solution of (1.8) for this case and also give
a simpler proof for Theorem 1.1. Our result extends (1.6) of
Jaschke (2003) to matrix-valued case. As an application, we reduce
some $\R^n$-valued inhomogeneous \emph{nonlinear} SDEs  to random
differential equations (RDEs) --- differential equations with
random coefficients. This facilitates pathwise study of solutions
and is an important step in the     context of random dynamical
systems approaches; see  Arnold (1998).

\section{Matrix-valued Inhomogeneous Linear SDEs }
\setcounter{equation}{0} \setcounter{Assumption}{0}
\setcounter{Theorem}{0} \setcounter{Proposition}{0}
\setcounter{Corollary}{0} \setcounter{Lemma}{0}
\setcounter{Definition}{0} \setcounter{Remark}{0}

L\'eandre (1985) obtained the following result about stochastic
equation (1.7). If we denote by $V$ the inverse of ${\cal E}(L)$,
then $V$ is the solution of the following equation:
$$ V(t)=I+\int_0^tV(s-) d W(s),  $$
where
$$W(t)=-L(t)+\langle L^c, L^c\rangle(t)+\sum_{0<s\le t}(1+\Delta
L(s))^{-1}(\Delta L(s))^2.$$ That means ${\cal E}(L)^{-1}={\cal
E}(W^\tau)^\tau$. We refer the reader to Karandikar (1991) for a
detailed proof of this result.

Now we will use this result of L\'eandre (1985) to solve the
following inhomogeneous stochastic integral equation
$$
X(t)=H(t)+\int_0^t (d L(s))X(s-), \eqno(2.1),
$$
where $L$ is a given $n\times n$ matrix-valued semimartingale with
$L 0=0$, such that $I+\Delta L(s)$ is invertible for all $s\in
[0,\infty]$, $H(t)$ is a given $n\times m$ matrix-valued adapted
cadlag process.

Our main result is the following.

\begin{Theorem} The unique solution of (2.1) is given by
$$X(t)=H(t)-{\cal E} (L)(t) \int_0^t (d {\cal E}
(L)(s)^{-1})H(s-).\eqno(2.2)$$ If $H$ is $n\times m$ matrix-valued
semimartingale, then $X(t)$ has the same expression as given by
(1.10) and (1.11).
\end{Theorem}

{\bf Proof}. \ \ We denote ${\cal E} (L)$ and ${\cal E} (L)^{-1}$ by
$U$ and $V$, respectively. We are going to show that the process
$$X(t)=H(t)-U(t)\int_0^t(dV(s))H(s-),$$
defined by (2.2), satisfies equation (2.1). Since $UV=I$, by the
integration by parts formula (see Karandikar (1991)) we get
\begin{eqnarray*}
0&=&d(U(t)V(t))=U(t-)dV(t)+(dU(t))V(t-)+d[U,V](t)\\
&=&d(U(t)V(t))=U(t-)dV(t)+dL(t)+d[U,V](t).
\end{eqnarray*}
Once again by the integration by parts formula, using the above
result  and the fact that $dU(t)=(dL(t))U(t-)$, we have
\begin{eqnarray*}d(X(t)-H(t))&=&-U(t-)dV(t)H(t-)-(dU(t))(\int_0^{t-} dV(s)
H(s-))-(d[U,V](t))H(t-)\\
&=&(dL(t))[H(t-)-U(t-)(\int_0^{t-} dV(s)H(s-))]=(dL(t))X(t-).
\end{eqnarray*}
This shows that the process $(X(t))$ defined by (2.2) satisfies
(2.1).

Now we assume that $(H(t))$ is an $n\times m$ matrix-valued
semimartingale. Using the notations in Section 1 we can verify that
$$G(t)=H(t)+[W,H](t). \eqno(2.3)$$
By the integration by parts formula, using (2.3) and the fact that
$dV(t)=V(t-)dW(t)$, we obtain
\begin{eqnarray*}
0&=&d(V(t)H(t))=V(t-)dH(t)+(dV(t))H(t-)+V(t-)d[W,H](t)\\
&=&V(t-)dG(t)+dL(t)+(dV(t))H(t-),
\end{eqnarray*}
from which we get
$$H(t)=U(t)\Big\{H(0)+\int_0^tV(s-)dG(s)+\int_0^tdV(s)H(s-)\Big\}.$$
Thus, if we let $(G(t))$ be defined by (2.3) then $(X(t))$ has the
expression of (2.2), and consequently it satisfies equation (2.1).
The proof of the theorem is complete. \hfill$\Box$

\section{An Application to Nonlinear SDEs}
\setcounter{equation}{0} \setcounter{Assumption}{0}
\setcounter{Theorem}{0} \setcounter{Proposition}{0}
\setcounter{Corollary}{0} \setcounter{Lemma}{0}
\setcounter{Definition}{0} \setcounter{Remark}{0}

In this section we will apply our results in Theorem 2.1 to reduce
an inhomogeneous \emph{nonlinear} SDE to a  RDE (random
differential equation).   Now we consider the following
$n$-dimensional nonlinear SDE (but with a linear multiplicative
noise term):
$$dX^i(t)=f^i(t,X(t))dt+\sum_{j=1}^nC^i_j(t)X^j(t)dB^j(t), \ \ X^i(0)=x^i,\eqno(3.1)$$
where $C(t)$ is an $n\times n$ matrix-valued measurable function,
$f(t,x)$ is a $\R^n$-valued measurable function on $[0,\infty)\times
\R^n$, and $B(t)=(B^1(t),\cdots, B^n(t))^\tau$ is a $n$-dimensional
Brownian motion. Put
$$ L^i_j(t)=\int_0^tC^i_j(s)dB^j(s), \ i.j=1,\cdots, n; \ \ H(t)=X(0)+\int_0^tf(s,X(s))ds.\eqno(3.2)$$
Then (3.1) can be rewritten in the form of  (2.1). For such $L$,
the equation (1.9) is reduced to the following linear equation:
$$dU^i_j(t)=\sum_{k=1}^nC^i_k(t)U^k_j(t)dB^k(t), \ \ U^i_j(0)=\delta^i_j, \eqno(3.3)$$
In the present case we have $G(t)=H(t)$. So according to (2.2),
the solution of (3.1) can be expressed as
$$X(t)=U(t)\Big \{ x +\int_0^t  U(t)^{-1}f(s,X(s))ds\Big \},\eqno(3.4)$$
where $U$ is the solution of (3.3). Unfortunately, even in this case
we are not able to give an explicit expression for $U(t)$. Let
$Y(t)=U(t)^{-1}X(t)$, then
$$Y(t)=\Big \{ x +\int_0^t U(s)^{-1}f(s,U(s)Y(s))ds\Big \}.\eqno(3.5)$$
This is the integral formulation of a RDE (random differential
equation).  Once we have sample path solution $Y(t)$ of this
transformed RDE , we obtain the solution of the original SDE via
$X(t)=U(t)Y(t)$.

%%%%%%%%%%%%%%%%%%%%%%%%%%%%%%%%%%%%%%%%%%%%%%%%%%%%%%%%%%%%%%%%%%%%%

\vspace{10mm} \noindent{\Large\bf References} {%\footnotesize
\begin{description}\baselineskip=5mm

\item{[1]}  Arnold, L. {\em Random Dynamical Systems.}
Springer-Verlag, New York, 1998.

\item{[2]} Dol\'eans-Dade, C. (1970): Quelques applications de la
formule de changement de variables pour les semimartingales, {Z.
Wahrsch. verw. Gebiete} 16, 181-194.

\item{[3]} Emery, M. (1978): Stablit\'e des solution des equations
differentielles stochastiques: application aux integrales
multiplicatives stochastique. {Z. Wahrsch. verw. Gebiete} 41,
241-262.

\item{[4]} Jacod, J. (1982): Equations diff\'erentielles
lin\'eares: la methode de variation des constantes, S\'eminaire de
Probabilit\'es XVI, LN in Math. 920, Springer-Verlag, 442-446.

\item{[5]} Jaschke, S. (2003): A note on the inhomogeneous linear
stochastic differential equation, {\it Insurance: Mathematics and
Finance} 32, 461-464.

\item{[6]} Karandikar, R.L. (1991): Multiplicative decomposition
of nonsingular matrix valued semimartingales, S\'eminaire de
Probabilit\'es XXV, LN in Math. 1485, Springer-Verlag, 262-269.

\item{[7]} L\'eandre, R. (1985) Flot d'une \'equation
diff\'erentielle stochastique, S\'eminaire de Probabilit\'es XIX,
LN in Math., Springer-Verlag, 271-274.

\item{[8]} Melnikov, A.V. and Shiryaev, A.N. (1996): Criteria for
the absence of arbitrage in the financial market, {\it Frontiers
in Pure and Appl. Probab. II}, Shiryaev, A.N. (Eds.), TVP Science
Publishers, Moscow, 121-134.

\item {[9]} Oksendal B. (1998): \emph{Stochastic differential
equations}, 5th Edition, Springer-Verlag.

\item{[10]} Protter, P. (2005): Stochastic integration and
differential equations, 2nd Edition, Springer-Verlag, New York.

\item{[11]} Revuz, D. and Yor, M. (1999): {\it Continuous
martingales and Brownian motion}, 3rd edition, Springer-Verlag,
Berlin.

\item{[12]} Yoeurp, C., Yor, M. (1977): Espace orthogonal\`a une
semimartingale, Unpublished.

\end{description}

\end{document}